\documentclass[reqno,a4paper]{amsart}

\usepackage{amsmath,amsthm,amsfonts,amssymb}
\usepackage{graphicx}
\usepackage{epstopdf}
\usepackage{enumerate}
\allowdisplaybreaks

\newtheorem{theorem}{Theorem}

\theoremstyle{definition}

\theoremstyle{remark}
\newtheorem{remark}{Remark}

\newcommand{\ct}{\tilde{c}}
\newcommand{\Ct}{\widetilde{C}}
\newcommand{\pt}{\tilde{p}}
\newcommand{\Pt}{\widetilde{P}}
\newcommand{\Dt}{\widetilde{D}}
\newcommand{\qt}{\tilde{q}}
\newcommand{\Qt}{\widetilde{Q}}
\newcommand{\at}{\tilde{a}}
\newcommand{\bt}{\tilde{b}}
\newcommand{\rt}{\tilde{r}}
\newcommand{\Rt}{\widetilde{R}}

\begin{document}
\date{}

\author{M.~Kazarian,~P.~Zograf}
\address{Steklov Mathematical Institute\\8 Gubkin St.\\Moscow 119991 Russia, and Department of Mathematics\\ National Research University Higher School of Economics\\Usacheva Str. 6\\ Moscow 119048 Russia}
\email{kazarian{\char'100}mccme.ru}
\address{St.Petersburg Department of the Steklov Mathematical Institute\\
Fontanka 27\\ St. Petersburg 191023 Russia, and
Chebyshev Laboratory of St. Petersburg State University\\
14th Line V.O. 29B\\
St.Petersburg 199178 Russia}
\email{zograf{\char'100}pdmi.ras.ru}

\title[Map and hypermap enumeration]{Rationality in map and hypermap enumeration by genus}


\begin{abstract}
Generating functions for a fixed genus map and hypermap enumeration become rational after a simple explicit change of variables. Their numerators are polynomials with integer coefficients that obey a differential recursion, and denominators are products of powers of explicit linear functions.
\end{abstract}

\maketitle

\section{Introduction}
By a {\em map} or a {\em ribbon graph} we understand a finite connected graph with prescribed cyclic orders of half-edges at each vertex. It also can be realized as the 1-skeleton of a polygonal partition of a closed orientable surface. The genus $g$ of a map (ribbon graph) satisfies the Euler formula
$$2-2g=\#v-\#e+\#f\,,$$
where $\#v,\#e,\#f$ are the numbers of vertices, edges and faces of the map respectively.
By a {\em hypermap} we understand a bicolored map, i.~e. a map whose faces are properly colored in two colors (say, white and black) so that no adjacent faces have the same color. The dual graph to a hypermap is a bipartite ribbon graph, or a {\em Grothendieck's ``dessin d'enfant''}
\footnote{As observed by Grothendieck, the absolute Galois group ${\rm Gal}(\overline{\mathbb{Q}}/\mathbb{Q})$ naturally acts on dessins (hypermaps); we refer the reader to \cite{LZ} for details.}.

We are interested in the weighted count of maps and hypermaps, where the weights are reciprocal to the orders of the corresponding automorphism groups. This is equivalent to counting {\em rooted} maps and hypermaps (i.~e. those with a marked half-edge). The passage from the rooted count to the unrooted one is known, cf. \cite{MN1}, \cite{MN2}.

Denote by $\ct_{g,n}$ (resp. $c_{g,n}$) the number of rooted maps (resp. hypermaps) of genus $g$ with $n$ edges (darts), and consider the genus $g$ generating functions
\begin{align}
&\Ct_g(s)=\sum_{n=2g}^\infty \ct_{g,n}s^n\,\quad g\geq 0\;,\\
&C_g(s)=\sum_{n=2g+1}^\infty c_{g,n}s^n\,\quad g\geq 0\;.
\end{align}
The classical problem that goes back to Tutte \cite{T} (or even earlier) is to compute the numbers $c_{g,n}$ and $\ct_{g,n}$. Effective algorithms for computing these numbers first appeared in \cite{WL} (for maps) and in \cite{Z} (for hypermaps)\footnote{The effective enumeration of 1-vertex maps was obtained in \cite{HZ}, and of 1-vertex hypermaps -- in \cite{J} and, independently, in \cite{A}. Enumeration of 1-vertex maps (or genus $g$ gluings of a $2n$-gon) was a crucial point in computing the Euler characteristic of the moduli space of algebraic curves in \cite{HZ}.}. Recursions for the numbers $c_{g,n}$ and $\ct_{g,n}$, and differential equations for the generating functions $C_g(s)$ and $\Ct_g(s)$ were first obtained in \cite{KZ} (cf. also \cite{CC} for an alternative approach to map enumeration).

In this note we show that the generating functions $C_g(s)$ and $\Ct_g(s)$ become rational functions after simple explicit changes of variable $s$. Their numerators are then polynomials with integer coefficients that obey a differential recursion, and denominators are products of powers of explicit linear functions.

\section{Main results}

We start with the case of hypermaps (Grothendieck's dessins d'enfants).

\begin{theorem}
Under the substitution $s=t(1-2t)$ we have
\begin{align}
&C_0(t(1-2t))=\frac{t(1-3t)}{(1-2t)^2}\,, \nonumber\\
&C_1(t(1-2t))=\frac{t^3}{(1-t)(1-4t)^2}\,, \nonumber\\
&C_g(t(1-2t))=\frac{P_g(t)}{(1-t)^{4g-3}(1-4t)^{5g-3}}\,,\quad g\geq 2,\label{hp}
\end{align}
where $P_g(t)=\sum_{i=2g+1}^{9g-7}p_{g,i}\,t^i$ is a polynomial with integer coefficients and $p_{g,2g+1}=\frac{(2g)!}{g+1}$. The polynomials $P_g(t)$ can be computed recursively by Eq.~(\ref{hint}).
\end{theorem}

\begin{remark}
The polynomials $P_g(t)$ for $g=2,3$ are:
\begin{align*}
P_2(t)&=8t^5-92t^6+464t^7-1316t^8+2204t^9-2048t^{10}+816t^{11},\\
P_3(t)&=180t^7-3648t^8+35424t^9-218944t^{10}+958160t^{11}-3102528t^{12}\\
      &+7503664t^{13}-13310768t^{14}+16365216t^{15}-11823680t^{16}+117916t^{17}\\
			&+6614784t^{18}-6008320t^{19}+1823744t^{20}.
\end{align*}
In principle, they can be computed for much larger values of $g$.
\end{remark}

\begin{remark}
Similar result was independently obtained in \cite{GW} by a different (more complicated) method.
\end{remark}

\begin{proof}
To prove the theorem, we recall a specialization of the Kadomtsev-Petviashvili (KP) equation for the hypermap count derived in \cite{KZ}:
\begin{align}
(sC_g)'&=3(2s^2C_g'+sC_g)+3s^3C_g'+s^3(s(sC_{g-1})')''\nonumber\\
&+s^3\sum_{i=0}^g(4C_i+6sC_i')C_{g-i}'+2s\delta_{g,0}\,,\label{h}
\end{align}
where the prime $'$ stands for the derivative $\frac{d}{ds}$.
This equation is just the differential form of the recursion (11) in \cite{KZ} for $t=u=v=1$:
\begin{align}
(n+1)c_{g,n}&=3(2n-1)\,c_{g,n-1}+(n-2)\,c_{g,n-2}+(n-1)^2(n-2)\,c_{g-1,n-2} \nonumber\\
            &+\sum_{i=0}^g\sum_{j=1}^{n-3}(4+6j)(n-2-j)\,c_{i,j}\,c_{g-i,n-2-j}\,.
\end{align}
For $g\geq 1$ we can further rewrite Eq.~(\ref{h}) as the differential recursion
\begin{align}
&(s-6s^2-3s^3-4s^3C_0-12s^4C_0')\,C_g'+(1-3s-4s^3C_0')C_g \nonumber\\
&=s^5C_{g-1}'''+5s^4C_{g-1}''+4s^3C_{g-1}'+s^3\sum_{i=1}^{g-1}(4C_i+6sC_i')C_{g-i}'\,.\label{dh}
\end{align}
For $g=0$ we get an ordinary differential equation that can be solved explicitly:
\begin{align*}
C_0(s)=\frac{-1+12s-24s^2+(1-8s)^{3/2}}{32s^2}\;.
\end{align*}
It is easy to see that the substitution $s=t(1-2t)$ considerably simplifies $C_0$ and makes it a rational function, namely
$$C_0(t(1-2t))=\frac{t(1-3t)}{(1-2t)^2}$$
(cf. \cite{W}). Substituting $s=t(1-2t)$ in (\ref{dh}), we get
\begin{align}
&t(1-t)^2(1-2t)\,\dot{C}_g+(1-t)(1-2t+4t^2)\,C_g \nonumber\\
&=t^3(1-2t)^3\left(D_tC_{g-1}+\frac{1}{(1-4t)^2}\,\sum_{i=1}^{g-1}\left(4(1-4t)C_i+6t(1-2t)\dot{C}_i\right)\dot{C}_{g-i}\right)\,,\label{dht}
\end{align}
where $\dot{C}_g=\frac{d}{dt}C_g(t(1-2t))$ and
\begin{align}
D_t=\frac{t^2(1-2t)^2}{(1-4t)^3}\cdot\frac{d^3}{dt^3}&+\frac{t(1-2t)(5-28t+56t^2)}{(1-4t)^4}\cdot\frac{d^2}{dt^2}\nonumber\\
&+\frac{4(1-11t+58t^2-144t^3+144t^4)}{(1-4t)^5}\cdot\frac{d}{dt}\,.
\end{align}
Assuming that $C_0,\ldots, C_{g-1}$ are known, we can think of (\ref{dht}) as an ODE for $C_g$. The integrating factor for this equation is
$\frac{1-t}{(1-2t)^3}$, so that we get from (\ref{dht})
\begin{align*}
\frac{d}{dt}&\left(\frac{t(1-t)^3}{(1-2t)^2}\,C_g\right)\nonumber\\
&=t^3(1-t)\left(D_tC_{g-1}+\frac{1}{(1-4t)^2}\,\sum_{i=1}^{g-1}\left(4(1-4t)C_i+6t(1-2t)\dot{C}_i\right)\dot{C}_{g-i}\right)\,,
\end{align*}
or, equivalently,
\begin{align}
&C_g(t(1-2t)=\frac{(1-2t)^2}{t(1-t)^3}\nonumber\\
&\times\int{t^3(1-t)\left(D_tC_{g-1}+\frac{1}{(1-4t)^2}\,\sum_{i=1}^{g-1}\left(4(1-4t)C_i+6t(1-2t)\dot{C}_i\right)\dot{C}_{g-i}\right)dt}\,.\label{hint}
\end{align}
Since by definition $C_g(0)=0$ for all $g\geq 0$, Eq.~(\ref{hint}) determines $C_g$ uniquely in terms of $C_0,\ldots, C_{g-1}$.
In particular, this equation immediately yields
\begin{align*}
&C_1(t(1-2t))=\frac{t^3}{(1-t)(1-4t)^2}\,,\\
&C_2(t(1-2t))=\frac{8t^5-92t^6+464t^7-1316t^8+2204t^9-2048t^{10}+816t^{11}}{(1-t)^5(1-4t)^7}\,.
\end{align*}
Let us show that $C_g(t(1-2t))$ has the form (\ref{hp}) for any $g\geq 3$. We will use the elementary formula
\begin{align}
\frac{d}{dt}\left(\frac{t^\alpha}{(1-t)^\beta (1-4t)^\gamma}\right)=\frac{\alpha t^{\alpha-1}+(-5\alpha+\beta+4\gamma)t^\alpha+4(\alpha-\beta-\gamma)t^{\alpha+1}}{(1-t)^{\beta+1} (1-4t)^{\gamma+1}}\,.\label{ef}
\end{align}
Then we have
\begin{align}
D_tC_{g-1}=\frac{(2g-1)(2g)^2p_{g-1,2g-1}t^{2g-2}+\ldots-256p_{g-1,9g-16}t^{9g-9}}{(1-t)^{4g-4} (1-4t)^{5g-2}}\label{hdo}
\end{align}
and
\begin{align}
\frac{1}{(1-4t)^2}\,&\sum_{i=1}^{g-1}\left(4(1-4t)C_i+6t(1-2t)\dot{C}_i\right)\dot{C}_{g-i}\nonumber\\
&=\frac{r_gt^{2g+1}+\ldots+256p_{g-1,9g-16}t^{9g-9}}{(1-t)^{4g-4} (1-4t)^{5g-2}}\,,\label{hsum}
\end{align}
where $r_g$ is some constant. Notice that the top degree term in the numerator on the right hand side of (\ref{hsum}) comes entirely from the product $C_1\dot{C}_{g-1}$.
Multiplying both sides of (\ref{hdo}) and (\ref{hsum}) by $t^3(1-t)$ and taking their sum we get that the integrand in (\ref{hint}) has the form
\begin{align}
\frac{Q_g(t)}{(1-t)^{4g-5} (1-4t)^{5g-2}}\,,
\end{align}
where $Q_g(t)=\sum_{i=2g+1}^{9g-7}q_{g,i}t^i$ is a polynomial with $q_{g,2g+1}=(2g-1)(2g)^2p_{g-1,2g-1}$. Therefore, we can rewrite (\ref{hint}) as
\begin{align}
C_g(t(1-2t))=\frac{(1-2t)^2}{t(1-t)^3}\int{\frac{Q_g(t)}{(1-t)^{4g-5} (1-4t)^{5g-2}}dt}\,.\label{hq}
\end{align}
To perform integration in (\ref{hq}) we decompose the integrand into the sum
\begin{align}
\frac{Q_g(t)}{(1-t)^{4g-5} (1-4t)^{5g-2}}=a+\sum_{i=2}^{4g-5}\frac{a_i}{(1-t)^i}+\sum_{j=2}^{5g-2}\frac{b_j}{(1-4t)^j}\,.\label{hel}
\end{align}
Note that no terms of the form $\frac{a_1}{1-t}$ or $\frac{b_1}{1-4t}$ can appear in the right hand side of (\ref{hel}) because the Taylor series expansion of the left hand side of (\ref{hq}) has integer coefficients\footnote{We owe this observation to F.~Petrov.}. Integrating we obtain
\begin{align}
\int{\frac{Q_g(t)}{(1-t)^{4g-5} (1-4t)^{5g-2}}dt}=at+b+\sum_{i=1}^{4g-6}\frac{a_{i+1}}{i}\cdot\frac{1}{(1-t)^i}+\sum_{j=1}^{5g-3}\frac{b_{j+1}}{4j}\cdot\frac{1}{(1-4t)^j}\,,\label{hrat}
\end{align}
where the condition $C_g(0)=0$ implies
\begin{align*}
b=-\sum_{i=1}^{4g-6}\frac{a_{i}}{i}-\sum_{j=2}^{5g-3}\frac{b_{j+1}}{4j}\,.
\end{align*}
Multiplying the right hand side of (\ref{hrat}) by $(1-t)^{4g-6}(1-4t)^{5g-3}$ we get a polynomial of the form $R_g(t)=\sum_{i=2g+2}^{9g-8}r_{g,i}t^i$. To complete the proof we put $P_g(t)=\frac{(1-2t)^2}{t}R_g(t)$ and notice that $p_{g,2g+1}=\frac{(2g-1)(2g)^2}{2g+2}\,p_{g-1,2g-1}$. Moreover, we see that $t=1/2$ is a root of $P_g(t)$ of multiplicity 2 provided $g\geq 2$.\footnote{Numerically, we also have $p_{g,9g-7}\neq 0,\; P_g(1)\neq 0,\; P_g(1/4)\neq 0$. In principle, this can be verified along the same lines as above, but computations become too cumbersome to reproduce them here.}
\end{proof}

Now we continue with map enumeration.

\begin{theorem} Under the substitution $s=t(1-3t)$ we have
\begin{align}
&\Ct_0(t(1-3t))=\frac{1-4t}{(1-3t)^2}\,, \nonumber\\
&\Ct_1(t(1-3t))=\frac{t^2}{(1-2t)(1-6t)^2}\,, \nonumber\\
&\Ct_g(t(1-3t))=\frac{\Pt_g(t)}{(1-2t)^{3g-2}(1-6t)^{5g-3}}\,,\quad g\geq 2,\label{mp}
\end{align}
where $\Pt_g(t)=\sum_{i=2g}^{8g-6}\pt_{g,i}\,t^i$ with $\pt_{g,2g}=\frac{(4g-1)!!}{2g+1}$.
The polynomials $\Pt_g(t)$ can be computed recursively by Eq.~(\ref{mint}).
\end{theorem}

\begin{remark}
The polynomials $\Pt_g(t)$ for $g=2,3$ are:
\begin{align*}
\Pt_2(t)&=21t^4-336t^5+2334t^6-9108t^7+21177t^8-27756t^9+15876t^{10}, \\
\Pt_3(t)&=1485t^6-41184t^7+539073t^8-4483458t^9+26893989t^{10}-124232004t^{11}\\
        &+453861279t^{12}-1307353122t^{13}+2897271774t^{14}-4737605112t^{15}\\
				&+5355443952t^{16}-3723895296t^{17}+1197496224t^{18}.
\end{align*}
Like in the case of hypermaps, they can be computed for much larger values of $g$.
\end{remark}

\begin{proof}
The proof of Theorem 2 is very similar to that of Theorem 1, so we will outline only its main steps.
We recall a specialization of the Kadomtsev-Petviashvili (KP) equation for the map count derived in \cite{KZ}:
\begin{align}
(s\Ct_g)'&=4(2s^2\Ct_g'+s\Ct_g)+2s^3(2s(s\Ct_{g-1})'+s\Ct)''+s^2(2s(s\Ct_{g-1})'+s\Ct)'\nonumber\\
&+3s^2\sum_{i=0}^g(\Ct_i+2s\Ct_i')(\Ct_{g-i}+2s\Ct_{g-i}')+\delta_{g,0}\,,\label{m}
\end{align}
where the prime $'$ stands for the derivative $\frac{d}{ds}$.
This equation is just a differential form of the recursion (16) in \cite{KZ} for $t=u=1$:
\begin{align}
(n+1)\ct_{g,n}&=4(2n-1)\,\ct_{g,n-1}+(2n-1)(2n-3)(n-1)\,\ct_{g-1,n-2} \nonumber\\
            &+3\sum_{i=0}^g\sum_{j=0}^{n-2}(2j+1)(2(n-2-j)+1)\,\ct_{i,j}\,\ct_{g-i,n-2-j}\,.
\end{align}
For $g\geq 1$ we can further rewrite Eq.~(\ref{m}) as the differential recursion
\begin{align}
&(s-8s^2-12s^3\Ct_0-24s^4\Ct_0')\,\Ct_g'+(1-4s-6s^2\Ct_0-12s^3\Ct_0')\,\Ct_g \nonumber\\
&=4s^5\Ct_{g-1}'''+24s^4\Ct_{g-1}''+27s^3\Ct_{g-1}'+3s^2\Ct_{g-1}+3s^2\sum_{i=0}^g(\Ct_i+2s\Ct_i')(\Ct_{g-i}+2s\Ct_{g-i}')\,.\label{dm}
\end{align}
For $g=0$ we get an ordinary differential equation that can be solved explicitly:
\begin{align*}
\Ct_0(s)=\frac{-1+18s+(1-12s)^{3/2}}{54s^2}\;.
\end{align*}
It is easy to see that the substitution $s=t(1-3t)$ considerably simplifies $\Ct_0$ and makes it a rational function, namely
$$\Ct_0(t(1-3t))=\frac{1-4t}{(1-3t)^2}$$
(cf. \cite{T}). Substituting $s=t(1-3t)$ in (\ref{dm}), we get
\begin{align}
&t(1-2t)(1-3t)\,\dot{C}_g+(1-4t+6t^2)\Ct_g=t^2(1-3t)^2\nonumber\\
&\times\left(\Dt_t\Ct_{g-1}+3\,\sum_{i=1}^{g-1}\left(\Ct_i+\frac{t(1-3t)}{1-6t}\dot{\Ct}_i\right)\left(\Ct_{g-i}+\frac{t(1-3t)}{1-6t}\dot{\Ct}_{g-i}\right)\right)\,,\label{dmt}
\end{align}
where $\dot{\Ct}_g=\frac{d}{dt}\Ct_g(t(1-3t))$ and
\begin{align}
\Dt_t=\frac{4t^3(1-3t)^3}{(1-6t)^3}&\cdot\frac{d^3}{dt^3}+\frac{24t^2(1-3t)^2(1-9t+27t^2)}{(1-6t)^4}\cdot\frac{d^2}{dt^2}\nonumber\\
&+\frac{9t(1-3t)(3-56t+456t^2-1728t^3+2592t^4)}{(1-6t)^5}\cdot\frac{d}{dt}+3\,.
\end{align}
Assuming that $\Ct_0,\ldots, \Ct_{g-1}$ are known, we can think of (\ref{dmt}) as an ODE for $\Ct_g$. The integrating factor for this equation is
$\frac{t(1-2t)}{1-3t}$, so that we get from (\ref{dmt})
\begin{align}
&\Ct_g(t(1-3t)=\frac{1-3t}{t(1-2t)}\nonumber\\
&\times\int{t^2\left(\Dt_t\Ct_{g-1}+
3\,\sum_{i=1}^{g-1}\left(\Ct_i+\frac{t(1-3t)}{1-6t}\dot{\Ct}_i\right)\left(\Ct_{g-i}+\frac{t(1-3t)}{1-6t}\dot{\Ct}_{g-i}\right)\right)dt}\,.\label{mint}
\end{align}
Since by definition $\Ct_g(0)=0$ for all $g\geq 1$, Eq.~(\ref{mint}) determines $\Ct_g$ uniquely in terms of $\Ct_0,\ldots, \Ct_{g-1}$.
In particular, this equation immediately yields
\begin{align*}
\Ct_1(t(1-3t))=\frac{t^2}{(1-2t)(1-6t)^2}\,.
\end{align*}
Let us show that $\Ct_g(t(1-3t))$ has the form (\ref{mp}) for any $g\geq 2$.
Using an analogue of (\ref{ef}) we obtain, after some cancellations, that the integrand in (\ref{mint}) has the form
\begin{align}
\frac{\Qt_g(t)}{(1-2t)^{3g-2} (1-6t)^{5g-2}}\,,
\end{align}
where $\Qt_g(t)=\sum_{i=2g}^{8g-6}\qt_{g,i}t^i$ is a polynomial with $\qt_{g,2g}=(2g-1)(4g-1)(4g-3)\,p_{g-1,2g-2}$.
It can be further decomposed into the sum
\begin{align}
\frac{\Qt_g(t)}{(1-2t)^{3g-2} (1-6t)^{5g-2}}=\sum_{i=2}^{3g-2}\frac{\at_i}{(1-2t)^i}+\sum_{j=2}^{5g-2}\frac{\bt_j}{(1-6t)^j}\,.\label{mel}
\end{align}
Integrating it, we obtain
\begin{align}
\int{\frac{\Qt_g(t)}{(1-2t)^{3g-2} (1-6t)^{5g-2}}dt}=\bt+\sum_{i=1}^{3g-3}\frac{\at_{i+1}}{2i}\cdot\frac{1}{(1-2t)^i}+\sum_{j=2}^{5g-3}\frac{\bt_{j+1}}{6j}\cdot\frac{1}{(1-6t)^j}\,,\label{mrat}
\end{align}
where the condition $\Ct_g(0)=0$ implies
\begin{align*}
\bt=-\sum_{i=1}^{3g-3}\frac{\at_{i+1}}{2i}-\sum_{j=1}^{5g-3}\frac{\bt_{j+1}}{6j}\,.
\end{align*}
Multiplying the right hand side of (\ref{mrat}) by $(1-2t)^{3g-3}(1-6t)^{5g-3}$ we get a polynomial of the form $\Rt_g(t)=\sum_{i=2g+1}^{8g-6}\rt_{g,i}t^i$. To complete the proof we put $\Pt_g(t)=\frac{1-3t}{t}\Rt_g(t)$ and notice that $\pt_{g,2g}=\frac{(2g-1)(4g-1)(4g-3)}{2g+1}\,\pt_{g-1,2g-2}$. Moreover, we see that $t=1/3$ is a root of $\Pt_g(t)$.\footnote{Numerically, we also have $\pt_{g,8g-6}\neq 0,\; \Pt_g(1/2)\neq 0,\; \Pt_g(1/6)\neq 0$.}

\end{proof}

{\bf Acknowledgements.} This work was supported by the Russian Science Foundation: Theorem~1 by the grant 14-21-00035, and Theorem~2 by the grant 16-11-10316. We are grateful to F.~Petrov for a useful suggestion. We thank A.~Giorgetti, A.~Mednykh, R.~Nedela and T.~R.~S.~Walsh for discussions.


\begin{thebibliography}{00}

\bibitem{A} Adrianov, N.: An analog of the Harer--Zagier formula for unicellular bicolored maps, Func. Anal. Appl. {\bf 31}:3, 1--9 (1997).
\bibitem{CC} Carrell, S.R., Chapuy, G.: Simple recurrence formulas to count maps on orientable surfaces. 	
J.~Combinatorial Theory, Series A, {\bf 133}, 58--75 (2015).
\bibitem{GW} Giorgetti, A., Walsh, T.~R.~S.: Enumeration of hypermaps of a given genus, arXiv:1510.09019 (2015).
\bibitem{HZ} Harer, J., Zagier, D.: The Euler characteristic of the moduli space of curves. Invent. Math.
{\bf 85}:3, 457--485 (1986).
\bibitem{J} Jackson, D.M.: Counting cycles in permutations by group characters, with an application to a topological problem.
Trans. Amer. Math. Soc. {\bf 299}, 785--801 (1987).
\bibitem{KZ} Kazarian, M.,  Zograf, P.: Virasoro constraints and topological recursion for Grothendieck’s dessin counting.
Lett. Math. Phys. {\bf 105} (8), 1057-1084 (2015).
\bibitem{LZ} Lando, S. K., Zvonkin, A. K.: Graphs on surfaces and their applications.
Encyclopaedia of Mathematical Sciences {\bf 141}, Springer-Verlag, Berlin (2004).
\bibitem{MN1} Mednykh, A., Nedela, R.: Enumeration of unrooted maps of a given genus. J.~Combinatorial Theory B {\bf 96}:5, 706--729 (2006).
\bibitem{MN2} Mednykh, A., Nedela, R.: Enumeration of unrooted hypermaps of a given genus. Discrete Mathematics {\bf 310}:3,
518--526  (2010).
\bibitem{T} Tutte, W.~T.: A census of planar maps. Canadian Math.~J. {\bf 15}:2, 249-271 (1963).
\bibitem{W} Walsh, T.~R.~S.: Hypermaps versus bipartite maps. J.~Combinatorial Theory B {\bf 18}:2, 155--163 (1975).
\bibitem{WL} Walsh, T.~R.~S., Lehman, A.~B.: Counting rooted maps by genus. I. J.~Combinatorial Theory B {\bf 13}, 192--218 (1972).
\bibitem{Z} Zograf, P.: Enumeration of Grothendieck's dessins and KP hierarchy, Int. Math. Res. Notices {\bf 2015}, 13533-13544 (2015).
\end{thebibliography}
\end{document}